\documentclass[conference]{IEEEtran}
\usepackage{graphicx}
\usepackage{url}
\usepackage{amsmath}

\begin{document}

\title{Statistics of voltage drop in radial distribution circuits: a dynamic programming approach}


\author{\authorblockN{Konstantin S. Turitsyn}
\authorblockA{T-4 \& CNLS, Theoretical Division, Los Alamos National Laboratory, Los Alamos, NM, USA \\}
\authorblockA{Landau Institute for Theoretical Physics, Moscow, Russia}
}

\maketitle

\begin{abstract}
We analyze a power distribution line with high penetration of distributed generation and strong variations of power consumption and generation levels. In the presence of uncertainty the statistical description of the system is required to assess the risks of power outages. In order to find the probability of exceeding the constraints for voltage levels we introduce the probability distribution of maximal voltage drop and propose an algorithm for finding this distribution. The algorithm is based on the assumption of random but statistically independent distribution of loads on buses. Linear complexity in the number of buses is achieved through the dynamic programming technique. We illustrate the performance of the algorithm by analyzing a simple 4-bus system with high variations of load levels.
\end{abstract}

\IEEEpeerreviewmaketitle

\section{Introduction}

Ensuring the quality of delivered power is one of the main challenges faced by utility companies. Rapid advent of renewable generation and electric vehicle technologies will inevitably result in a significant increase in the variations of  power consumption and/or generation Keeping the voltage level within the industry constraints will become an even more formidable task for the utilities. One of the most challenging aspect of this problem is the uncertainty about the load structure in the feeder line. Lack of information about the loads makes the problem essentially probabilistic and thus requires more sophisticated techniques for analysis and control.

Traditional methods of controling the voltage level on distribution feeder lines include but are not limited to line regulators and capacitor banks \cite{06NW}. Whereas the line regulators based transformer tap changes are employed for controlling the voltage level in the beginning of the line, the switched capacitors are usually distributed in the middle of the line, and can be used for smoothing out the voltage drop curve via the reactive power injections \cite{89BWa,89BWc}. Penetration of intermittent renewable generators and increasing demand in power supply will require novel approaches for controlling the voltage level. One of the most promising ideas that has been proposed recently in this field is the distributed control of reactive power flows via the local inverters attached to renewable generators \cite{08BMOCKS,08LB,08WNRN}.

In a recent work \cite{09TSBC} we have shown that the simple control techniques of distributed inverters can be efficient in reducing the losses in radial distribution systems. However, the effect of inverters on the overall power quality remained an open question. One of the main obstacles to the proper assessment of the power quality in the system was the lack of theoretical methods of voltage drop analysis in the presence of high variations of loads. In this paper we attempt first steps of solving this problem by introducing a novel algorithm, that allows fast prediction of the probability distribution function of the maximal voltage drop in the feeder line. Specifically the algorithm allows one to transform the distribution of loads on individual buses into the distribution of maximal voltage drop.

The structure of this paper is the following: in the section \ref{sec:model} we introduce the power flow model based on classical DistFlow equations. Then, we formally introduce the main object this study: the maximal voltage drop. Extension of the model to include the stochastic variations of power consumption/generation are presented in the section \ref{sec:stat}. We finish by presenting the results of simulations and discussing the remaining challenges.

\section{Power flow model}\label{sec:model}
The power flows in the linear distribution line with $N$ buses can be described with the DistFlow recurrence equations \cite{89BWa,89BWc}:
\begin{eqnarray}\label{Peq}
 P_{k+1} = P_k - r_k \frac{P_k^2 + Q_k^2}{V_k^2} - p_{k+1} \\
 Q_{k+1} = Q_k - x_k \frac{P_k^2 + Q_k^2}{V_k^2} - q_{k+1} \label{Qeq}\\
 V_{k+1}^2  = V_k^2 - 2 (r_k P_k + x_k Q_k) + (r_k^2+x_k^2) \frac{P_k^2 + Q_k^2}{V_k^2}  \label{Veq0}
\end{eqnarray}
where $P_k,Q_k$ are real and reactive power flows from bus $k-1$ to $k$. $V_k$ is the voltage level on bus $k$, and $p_k,q_k$ are the values of real and reactive power consumption on bus $k$. $r_k+j x_k$ is the complex impendance of the link between the nodes $k$ and $k+1$. These equations have to be solved with two boundary conditions: fixed base voltage level at the beginning of the line $V_0$ and zero power flux  through the virtual link at the end of the line: $P_{N+1} = Q_{N+1} = 0$.

\begin{figure}
 \includegraphics[width=0.5\textwidth]{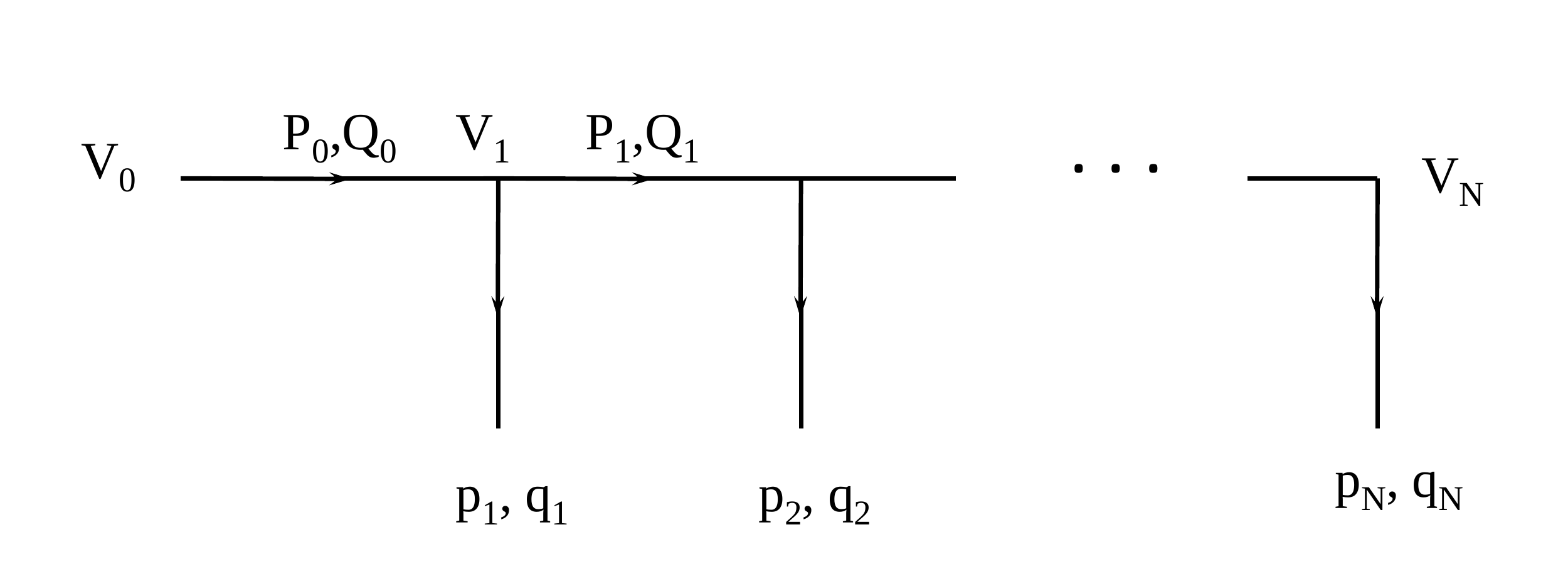}
\centering \caption{Schematic representation of the linear feeder line. Each bus is characterized by the voltage level $V_k$, consumption of real and reactive powers: $p_k,q_k$. The power flows between the buses are denoted by $P_k,Q_k$.
}\label{fig:scheme}
\end{figure}

For a class of low voltage distribution lines considered here the quadratic terms in the equations (\ref{Peq}-\ref{Veq0}) are negligibly small, so one can use the linearized power flow equations where quadratic terms are dropped. Moreover, for the sake of simplicity in this work we will restrict the analysis to the homogeneous networks where the ratio $x_k/r_k = \alpha$ is the same for all links. In this case one can combine the values of $P_k$ and $Q_k$ into a single variable $S_k = P_k + \alpha Q_k$ (and also introduce $s_k = p_k + \alpha q_k$). Moreover, as long as the variations of the voltage are small compared to its base value, one can use the approximation $V_k^2 = V_0^2 + 2V_0 (V_k-V_0)$. Using these approximations one can rewrite the DistFlow equations in the following simple form:
\begin{eqnarray}\label{Seq}
 S_{k} =  S_{k+1} + s_{k+1}\\
 V_k = V_{k+1} + \rho_k S_k \label{Veq}.
\end{eqnarray}
where we have introduced $\rho_k = r_k/V_0$.
\subsection{Maximal voltage drop}
The power quality in the line can be characterized by the maximal voltage drop level - the difference between the base voltage $V_0$ and the minimal voltage in all buses: $\Delta_0 = V_0 - \min_{k} V_k$. In order to find the value of $\Delta_0$ we introduce an intermediate object $\Delta_n$ characterizing the maximal voltage drop in the feeder line segment $n\dots N$:
\begin{equation}\label{DeltaEq}
 \Delta_n = V_n - \min_{n \leq k \leq N} V_k
\end{equation}
Note, that $\Delta_n$ is essentially non-negative quantity.  The recurrence equation for $\Delta_n$ can be easily derived from (\ref{Veq}) by noticing that there are two distinct scenarios. Whenever $V_{k+1}+\rho_k S_k > 0$ the voltage level on bus $k$ will be larger than the minimal one in the segment $k\dots N$ : $V_k =  V_{k+1} + \rho_k S_k > V_{k+1} -\Delta_{k+1} = \min_{k+1 \leq l \leq N} V_l$. Therefore, according to the definition (\ref{DeltaEq}) the value of $\Delta_k$  will be given by $\Delta_k = \Delta_{k+1} + \rho_k S_k$. In the second case when $\Delta_{k+1}+\rho_k S_k \leq 0$ the value of $V_k$ will be smaller than $\min_{k+1 \leq l \leq N} V_l$ and thus the minimal value of the voltage in the segment $k\dots N$ is achieved on bus $k$, and thus $\Delta_k = 0$. Formally one can express these two scenarios with the following equation:
\begin{equation}\label{Deq}
 \Delta_k =
\begin{cases} \Delta_{k+1} + \rho_k S_k & \text{if $\Delta_{k+1} + \rho_k S_k > 0$,} \\
0 &\text{if $\Delta_{k+1} + \rho_k S_k\le 0$.}
\end{cases}
\end{equation}
It is important to note, that although the equations (\ref{Seq},\ref{Veq}) for $S_k$ and $U_k$ are linear, the resulting equation (\ref{Deq}) for $\Delta_k$ is essentially non-linear because of the second case where $\Delta_k$ is reset to zero. This second case is realized only if the total power $S_k$ flowing through the link is negative. This can happen only in the presence of buses that inject power in the system. As long as $S_k$ variable incorporates both real and reactive power, the injection of $S$ flow can be associated either with distributed generators that inject real power or with capacitor banks that inject reactive power. In the absence of power injection in the system, the voltage is a monotonous function, and the problem of finding the maximal voltage drops reduces to calculation of $V_0 -V_N$ which is a linear function of the loads $s_k$. This problem can be solved analytically without any sophisticated algorithms, proposed in this paper. Here we focus on the nontrivial situation, and assume that there are buses that inject power in the line.

By solving the recurrence relations (\ref{Seq}) and (\ref{Deq}) backwards in $k$ with the intial conditions $S_{N+1} = 0, \Delta_{N+1} = 0$ one can easily find the value of maximal voltage drop $\Delta_0$ in the whole network. These equations form a basis for the statistical analysis of voltage drop in presence of uncertainties about load levels in the system.

\section{Statistical description of power flow}\label{sec:stat}
Whenever the precise values of $p_k,q_k$ are not known, it is not possible to find the precise value of the maximal voltage drop $\Delta_0$. Instead one has to develop a statistical approach to the problem. The traditional way of characterizing the uncertainties in the system is to study the probability density functions (PDF) of system state. Assuming that the loads on different buses change independently one can define the probability of observing $p_k + \alpha q_k = s$ via the PDF $\pi_k(s)$ \footnote{Throughout the paper we use the term probability distribution to refer to the probability density function: by definition the probability of observing $s_k \in [s-\delta/2,s+\delta/2]$ is equal to $\pi_k(s)\delta$ in the limit of $\delta\to 0$}. The origin of these distributions is not important for our analysis: they could be either prior ``Bayesian'' distributions based on the measurement history, or distributions derived from some statistical model of the load on bus $k$. The statistical independce assumption implies that the joint probability of observing $p_1+\alpha q_1 = s_1',\dots p_N + \alpha q_N = s_N'$ is given by the product $\pi(s_1'\dots s_N') = \prod_k \pi_k(s_k')$.

The goal of the statistical analysis is to transform the PDF of loads into the PDF of the maximal voltage drop $\Delta_0$ i.e. into the probability of observing the given value of maximal voltage drop. One can easily write the formal integral expression for the probability of having $\Delta_0 < \Delta$. This can be done by noting that $\Delta_0 < \Delta$ if and only if $V_k > V_0 - \Delta$ for all $0\leq k \leq N$. Given the formal solution $V_k = V_k(s_1\dots s_N)$ of (\ref{Seq},\ref{Veq}) one can formally write:
\begin{equation}
 \mathrm{Prob}(\Delta_0 < \Delta ) = \int \prod_k ds_k \pi(s_k) \theta(V_k(s_1\dots s_N)- V_0 + \Delta)
\end{equation}
where $\theta(x)$ is the Heaviside (Unit step) function. Here and throughout the text we assume that intergration is taken over the domain $(-\infty,\infty)$, unless the integration domain is specified explicitly. On practical level, this solution is unusable for large systems with $N\gg 1$, as long as there are no fast ``black-box" algorithms to evaluate this multidimensional integral in case $N \gg 1$. Even in the simplest case of uniform distributions $\pi_k(s_k)$ the problem is reduced to calculation of the volume of a highly dimensional polytope, for which state of the art algorithms require at least $O(N^4)$ operations \cite{06LV}.

It is therefore important to develop alternative approaches to computing the PDF of $\Delta_0$ that would require only linear in $N$ number of operations and would be applicable even for large systems. We propose a specific  algorithm of the kind that exploits the radial structure of the network and statistical independence of load levels on different buses. The main idea behind the algorithm is to transform the recurrence relations (\ref{Seq},\ref{Deq}) into the corresponding relations for the joint probability functions of $\Delta_k, S_k$ and solve them in an iterative way. This general strategy of reusing the previous computations via decomposing the problem in recurrence type relations is usually referred as dynamic programming \cite{01CLRS}.

Existence of recurrence relations (\ref{Seq},\ref{Deq}) for the pair of variables $S_k,\Delta_k$ which expresses the PDF $\Pi_k(S,\Delta)$ of the values of $S_k,\Delta_k$ via the PDF $\Pi_{k+1}(S,\Delta)$ of the values of $S_{k+1},\Delta_{k+1}$. This Chapman-Kolmogorov type relation can be written in a most compact way with the use of Dirac Delta Functions $\delta(x)$:
\begin{eqnarray}\label{PiEq}
 \Pi_k(S,\Delta) = \int dS' d\Delta' \pi_{k+1}(S-S')\Pi_{k+1}(S',\Delta')\times \nonumber\\
\left[\theta(\Delta'+\rho_{k}S)\delta(\Delta-\Delta' - \rho_k S) + \theta(-\Delta'-\rho_{k}S)\delta(\Delta)\right] 
\end{eqnarray}
where we have used $S',\Delta'$ to denote the values of $S_{k+1},\Delta_{k+1}$ respectively. This recurrence relation has to be solved with the initial condition $\Pi_{N+1}(S,\Delta) = \delta(S)\delta(\Delta)$. Formally this equation completes the construction of the algorithm, as one can iterate it $N$ times to obtain the PDF of maximal voltage drop in the system $\Delta_0$. For a given target precision, the complexity of performing single iteration depends only on the numerical discretization of $\Pi_k(S,\Delta)$  and does not depend on the total number of buses in the system. Therefore the total number of iterations grows like $O(N)$. However, implementation of a single iteration (\ref{PiEq}) can be obstructed by a non-analytic nature of the resulting joint PDF $\Pi_k(S,\Delta)$. Performing the first few iterations analytically, one can show that for continous distributions $\pi_k(s)$ the general expression for $\Pi_k(S,\Delta)$ can be decomposed in three parts:
\begin{equation}
  \Pi_k(S,\Delta) = \Pi_k^{c}(S,\Delta) + \Pi_k^{(1)}(S)\delta(\Delta) + \Pi_k^{(2)}(S)\delta(\Delta - 2 r_k S)
\end{equation}
where $\Pi_k^c(S,\Delta)$ is analytic part of the PDF, and $\Pi_k^{(1,2)}$ are the prefactors in front of the non-analytic ones. Substituting this decomposition in (\ref{PiEq}) one obtains the coupled system of equations for $\Pi_k^c(S,\Delta)$ and $\Pi_k^{(1,2)}$:
\begin{eqnarray}
 \Pi_k^c(S,\Delta) = \theta(\Delta)\int dS' \pi_{k+1}(S-S')\Pi_{k+1}^c(S',\Delta-2 r_k S) + \nonumber \\
 \frac{\theta(\Delta)\theta(\Delta-S)}{2 r_{k+1}}\pi_{k+1}\left(S-\frac{\Delta-2 r_k S}{2 r_{k+1}}\right)\Pi_{k+1}^{(2)}\left(\frac{\Delta-2 r_k S}{2 r_{k+1}}\right)\label{Pceq}\nonumber
\end{eqnarray}
\begin{eqnarray}
& \Pi_k^{(1)}(S) = \theta(-S)\int dS'\pi_{k+1}(S-S')\Pi_{k+1}^{(1)}(S') +  \nonumber \\
&\theta(-S)\int dS'\pi_{k+1}(S-S')\int_{-\infty}^{-2 r_k S} d\Delta' \Pi_{k+1}^c(S',\Delta') + \nonumber\\
&\theta(-S)\int_0^{- r_k S/r_{k+1}} dS'\pi_{k+1}(S-S')\theta(S')\Pi_{k+1}^{(2)}(S')\label{P1eq}\nonumber
\end{eqnarray}
\begin{eqnarray}\label{P2eq}
  \Pi_k^{(2)}(S) =  \theta(S)\int dS' \pi_{k+1}(S-S')\Pi_{k+1}^{(1)}(S')
\end{eqnarray}
This system of equations has to be solved backwards in $n$ with the initial condition $\Pi_{N+1}^c(S,\Delta) = 0,\Pi_{N+1}^{(1)}(S) = \delta(S)$ and $\Pi_{N+1}^{(2)}(S) = 0$. In the next section we will describe the results of solving this system for a sample $4$ bus system.
\section{Results}

\begin{figure}
\includegraphics[width=0.5\textwidth]{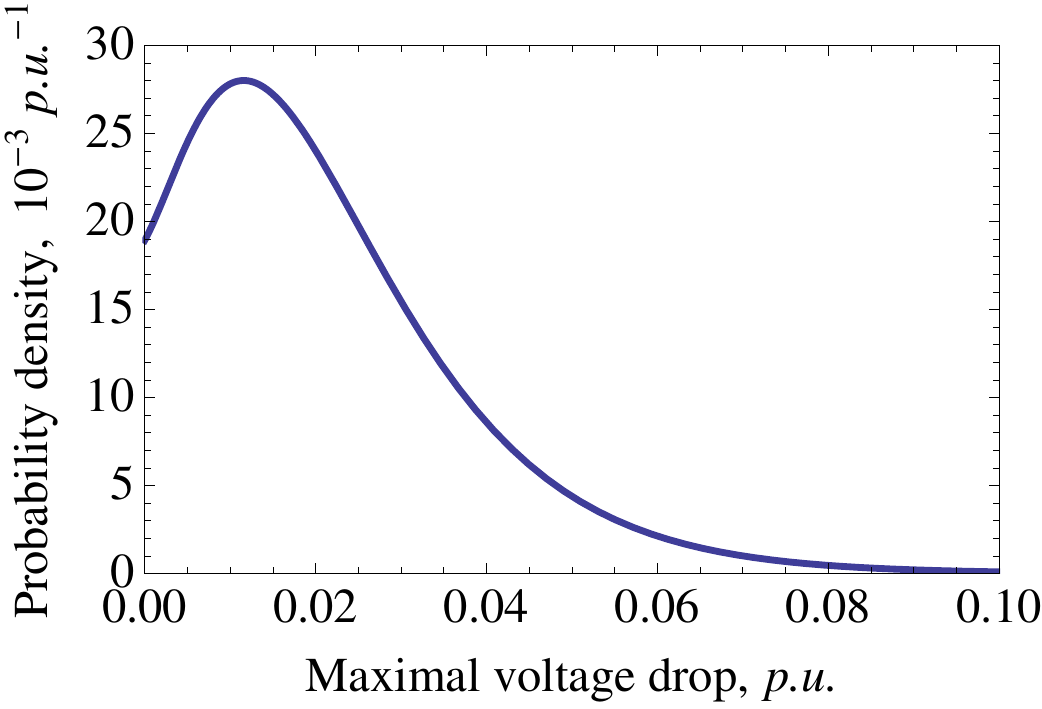}
\centering \caption{Probability distribution function of the maximal voltage drop $\Delta_0$.
}\label{fig:prob}
\end{figure}

Efficient implementation of the proposed algorithm suitable for large heterogeneous systems is a challenging task that requires a thoughtful selection of suitable discretization technique for the PDFs, that would allow fast evaluation of convolution integrals in (\ref{P2eq}) and adaptation of the discretization domain to the width of the distribution functions that will grow with each iteration. Analysis of possible approaches to these problems is beyond the scope of this paper. Instead, here we report our results of a ``proof of concept'' study of toy feeder line system consisting of $N=4$ buses where for the sake of simplicity we have set all $\rho_k$ to $\rho_k = 10^{-3} p.u./kW$ and the load distribution functions to
\begin{equation}
 \pi_k(s)  =  0.25\begin{cases} \exp(-s/3.0) & \text{if $s > 0$,} \\
\exp(1.0 s) &\text{if $s \le 0$.}
\end{cases}
\end{equation}
Note, that this distribution function implies that flux of real and reactive power can have an arbitrary sign. Although the probability of power consumption ($s>0$) is higher, there is also a finite probability of a given bus injecting the power in the line. The average value of power consumption was set to $\bar s = 2.0 kW$  whereas the standard deviation of power consumption is equal to $\sigma_s \approx 3.16 kW$, so the fluctuations of power consumption/generation are indeed strong in the system.

In order to avoid technical difficulties associated with discretization of the probability distribution functions, we have solved the equations (\ref{P2eq}) analytically using the Wolfram Mathematica computer algebra system. On each step the functions $\Pi_k^c(S,\Delta),\Pi_k^{(1,2)}(S)$ were represented as a piecewise analytical functions, and the convolutions (\ref{P2eq}) were performed rigorously without any numerical approximations. Although the number of terms in the analytical expressions was growing quite rapidly, it was possible to find the expressions for $\Pi_0^c(S,\Delta)$ and $\Pi_0^{(1,2)}(S)$ for a system with $N=4$ buses. These expressions were used to calculate several statistical characteristics of the system.

The probability of exceeding the voltage level constraints can be calculated from the analytical part of the probability distribution of the real voltage drop $\Delta_0$ or equivalently the probability density in the region of positive voltage drops:
\begin{equation}
 P(\Delta_0| \Delta_0>0) =\Pi_0^{(2)}(\Delta_0) + \int dS' \Pi_0^c(S', \Delta_0)
\end{equation}
This corresponding result is presented on the figure \ref{fig:prob}. As one can see the distribution is centered around characteristic value of $\bar \Delta_0\approx 4.1 p.u.$, however the variations of the typical voltage drop are very high, for instance the probability of exceeding twice the average value $\Delta_0 > 2 \bar \Delta_0$ is about $13\%$. This result confirms the central thesis of the study: variations in load levels result in large variations of the voltage drop along the line and thus significantly increase the probability of exceeding the limits set by regulators.

\begin{figure}
\includegraphics[width=0.5\textwidth]{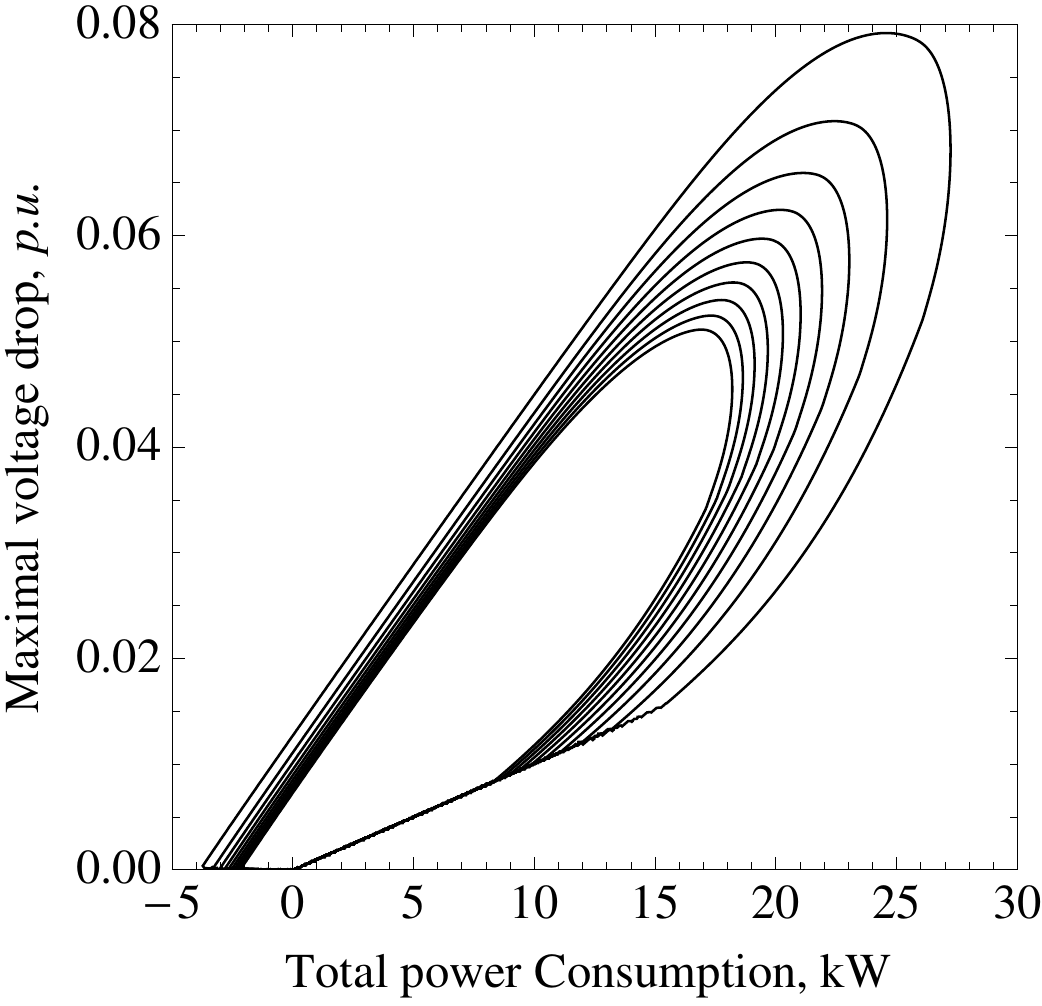}
\centering \caption{Contour plot of the joint probability distribution $\Pi_0^c(S,\Delta)$ of the total load $S_0$ and maximal voltage drop $\Delta_0$.
} \label{fig:contour}
\end{figure}

Another interesting object is the joint distribution of the total power consumption in the system $S_0$ and  the maximal voltage drop $\Delta_0$. The corresponding distribution  $\Pi_0^c(S_0,\Delta_0)$ is shown on the figure \ref{fig:contour}. As one can see there is a strong correlation between $S_0$ and $\Delta_0$. This correlation can be potentially used for designing the voltage control techniques that respond to the total load $S_0$ and adjust the voltage $V_0$ to suppress the risk of exceeding the allowev voltage drop limits. Statistical analysis of various control techniques is an interesting problem that can be approached with the techniques introduced in this article.

\section{Discussion}
Algorithm presented in this paper is a first step on a long path of development a new generation of tools that would allow the power system designers and controllers to assess the risks associated with fluctuations of load levels and uncertainties in the system. Fast ways of calculating the probabilities of high voltage drops can be utilized in system design optimization packages or stochastic control processes. One of the specific applications of this algorithm could be the analysis of the effect of distributed renewable generation on the power quality in the distribution system. The probability distribution functions of voltage drop can be used for determination of the critical penetration levels of renewable generators that can be sustained by the current lines. It can be also used for comparison of different techniques of voltage control with reactive power generators: capacitor banks \cite{89BWa} or small-scale inverters attached to photovoltaic elements \cite{09TSBC}.

However, there are still several issues that have to be dealt with before the algorithm can be put to use. The most important of them is the choice of suitable discretization technique, that would allow efficient evaluation of  integrals (\ref{P2eq}). The main obstacle here is the piecewise-analytic structure of the convolution kernels that will result in piecewise continous structure of the distribution functions. This structure suggests that the probability distributions should be discretized with finite element partitioning of the domain. However, the composition of domains will have to be adaptable, as the structure of the distribution function will change with each iteration. For practical problems it is also important to extend the current algorithm to analyze the joint statistics not only of the minimal but also the maximal levels of voltage along the line. Although on formal mathematical level this extension is straightforward, the numerical realization will become more complicated, as the distribution functions will depend on three variables: power flow $S$ and maximal and minimal voltage drops.

\section{Conclusion}
We have presented a dynamic programming algorithm for calculating the probability distribution function of the maximal voltage drop in the linear feeder line. The complexity of the algorithm scales linear with the number of buses in the system, which makes it superior in comparison to the standard black-box type algorithm for solving similar problems. The proposed algorithm requires iterative application of integral operators that relate the probability distributions of local voltage drops on neighboring nodes. Efficient implementation of this convolution is a difficult problem that will be hopefully solved in forthcoming works. In this paper we have tested the algorithm on a toy $4$-bus homogeneous system. The resulting distribution functions show that the variations in load lead to variations in voltage drop. Moreover they also show strong correlation between the voltage drop and total load of the system.

\section*{Acknowledgment}
We are indebted to Scott Backhaus for drawing the attention to the problem and thank Misha Chertkov and Petr Sulc for useful discussions. This work has been supported by the Oppenheimer fellowship from Los Alamos National Laboratory.

\bibliographystyle{IEEEtran}

\bibliography{feederSIBIRCON2010}

\end{document}